\documentclass[12pt]{article}%
\usepackage{graphicx}
\usepackage[intlimits]{amsmath}
\usepackage{latexsym}
\usepackage{amsfonts}
\usepackage{amssymb}%
\makeatletter
\newfont{\footsc}{cmcsc10 at 8truept}
\newfont{\footbf}{cmbx10 at 8truept}
\newfont{\footrm}{cmr10 at 10truept}
\newtheorem{theorem}{Theorem}

\newtheorem{corollary}[theorem]{Corollary}

\newtheorem{lemma}[theorem]{Lemma}

\newtheorem{proposition}[theorem]{Proposition}

\newenvironment{proof}[1][Proof]{\noindent{\textbf {#1}  }}  {\hfill$\Box$\bigskip}

\begin{document}

\title{Linear combinations of graph eigenvalues}
\author{Vladimir Nikiforov\\Department of Mathematical Sciences, University of Memphis, \\Memphis TN 38152, USA, email: \textit{vnkifrv@memphis.edu}}
\maketitle

\begin{abstract}
Let $\mu_{1}\left(  G\right)  \geq\ldots\geq\mu_{n}\left(  G\right)  $ be the
eigenvalues of the adjacency matrix of a graph $G$ of order $n,$ and
$\overline{G}$ be the complement of $G.$

Suppose $F\left(  G\right)  $ is a fixed linear combination of $\mu_{i}\left(
G\right)  ,$ $\mu_{n-i+1}\left(  G\right)  ,$ $\mu_{i}\left(  \overline
{G}\right)  ,$ and $\mu_{n-i+1}\left(  \overline{G}\right)  ,$ $1\leq i\leq
k.$ We show that the limit
\[
\lim_{n\rightarrow\infty}\frac{1}{n}\max\left\{  F\left(  G\right)  :v\left(
G\right)  =n\right\}
\]
always exists. Moreover, the statement remains true if the maximum is taken
over some restricted families like \textquotedblleft$K_{r}$%
-free\textquotedblright\ or \textquotedblleft$r$-partite\textquotedblright\ graphs.

We also show that
\[
\frac{29+\sqrt{329}}{42}n-25\leq\max_{v\left(  G\right)  =n}\mu_{1}\left(
G\right)  +\mu_{2}\left(  G\right)  \leq\frac{2}{\sqrt{3}}n,
\]
answering in the negative a question of Gernert.



\end{abstract}

\section{Introduction}

Our notation is standard (e.g., see \cite{Bol98}, \cite{CDS80}, and
\cite{HoJo88}); in particular, all graphs are defined on the vertex set
$\left[  n\right]  =\left\{  1,\ldots,n\right\}  $ and $\overline{G}$ stands
for the complement of $G.$ We order the eigenvalues of the adjacency matrix of
a graph $G$ of order $n$ as $\mu_{1}\left(  G\right)  \geq\ldots\geq\mu
_{n}\left(  G\right)  .$

Suppose $k>0$ is a fixed integer and $\alpha_{1},\ldots,\alpha_{k},\beta
_{1},\ldots,\beta_{k},\gamma_{1},\ldots,\gamma_{k},\delta_{1},\ldots
,\delta_{k},$ are fixed reals. For any graph $G$ of order at least $k,$ let%
\[
F\left(  G\right)  =\sum_{i=1}^{k}\alpha_{i}\mu_{i}\left(  G\right)
+\beta_{i}\mu_{n-i+1}\left(  G\right)  +\gamma_{i}\mu_{i}\left(  \overline
{G}\right)  +\delta_{i}\mu_{n-i+1}\left(  \overline{G}\right)  .
\]

For a given graph property $\mathcal{F}$, i.e., a family of graphs closed
under isomorphism, it is natural to look for $\max\left\{  F\left(  G\right)
:G\in\mathcal{F}\text{, }v\left(  G\right)  =n\right\}  .$ Questions of this
type have been studied, here is a partial list:%
\[%
\begin{array}
[c]{ll}%
\max\left\{  \mu_{1}\left(  G\right)  +\mu_{n}\left(  G\right)  :\text{
}G\text{ is }K_{r}\text{-free,\ }v\left(  G\right)  =n\right\}  & \text{Brandt
\cite{Bra98};}\\
\max\left\{  \mu_{1}\left(  G\right)  -\mu_{n}\left(  G\right)  :v\left(
G\right)  =n\right\}  & \text{Gregory, Hershkowitz, Kirkland \cite{GHK01};}\\
\max\left\{  \mu_{1}\left(  G\right)  +\mu_{2}\left(  G\right)  :v\left(
G\right)  =n\right\}  & \text{Gernert \cite{Ger};}\\
\max\left\{  \mu_{1}\left(  G\right)  +\mu_{1}\left(  \overline{G}\right)
:v\left(  G\right)  =n\right\}  & \text{Nosal \cite{Nos70}, Nikiforov
\cite{Nik02};}\\
\max\left\{  \mu_{i}\left(  G\right)  +\mu_{i}\left(  \overline{G}\right)
:v\left(  G\right)  =n\right\}  & \text{Nikiforov \cite{Nik06}.}%
\end{array}
\]

One of the few sensible questions in such a general setup is the following
one: does the limit
\[
\lim_{n\rightarrow\infty}\frac{1}{n}\max\left\{  F\left(  G\right)
:G\in\mathcal{F}\text{, }v\left(  G\right)  =n\right\}
\]
exist? We show that, under some mild conditions on $\mathcal{F},$ this is
always the case.

For any graph $G=\left(  V,E\right)  $ and integer $t\geq1,$ write $G^{\left(
t\right)  }$ for the graph obtained by replacing each vertex $u\in V$ by a set
$V_{u}$ of $t$ independent vertices and joining $x\in V_{u}$ to $y\in V_{v}$
if and only if $uv\in E.$

Call a graph property $\mathcal{F}$ \emph{multiplicative} if : \emph{(a)}
$\mathcal{F}$ is closed under adding isolated vertices; \emph{(b)}
$G\in\mathcal{F}$ implies $G^{\left(  t\right)  }\in\mathcal{F}$ for every
$t\geq1$. Note that \textquotedblleft$K_{r}$-free\textquotedblright,
\textquotedblleft$r$-partite\textquotedblright, and \textquotedblleft any
graph\textquotedblright\ are multiplicative properties.

\begin{theorem}
\label{th1}For any multiplicative property $\mathcal{F}$ the limit
\begin{equation}
c=\lim_{n\rightarrow\infty}\frac{1}{n}\max\left\{  F\left(  G\right)
:G\in\mathcal{F}\text{, }v\left(  G\right)  =n\right\}  \label{eq1}%
\end{equation}
exists. Moreover,
\[
c=\lim\sup\left\{  \frac{1}{\left\vert G\right\vert }F\left(  G\right)
:G\in\mathcal{F}\right\}  .
\]

\end{theorem}

Note that, since the $\alpha_{i}$'s$,\beta_{i}$'s$,\gamma_{i}$'s$,$ and
$\delta_{i}$'s may have any sign, Theorem \ref{th1} implies that
\[
\lim_{n\rightarrow\infty}\frac{1}{n}\min\left\{  F\left(  G\right)
:G\in\mathcal{F}\text{, }v\left(  G\right)  =n\right\}
\]
exists as well.

Gernert \cite{Ger} (see also Stevanovic \cite{Ste06}) has proved that the
inequality
\[
\mu_{1}\left(  G\right)  +\mu_{2}\left(  G\right)  \leq v\left(  G\right)
\]
holds if the graph $G$ has fewer than $10$ vertices or is one of the following
types: regular, triangle-free, thoroidal, or planar; he consequently asked
whether this inequality holds for any graph $G$. We answer this question in
the negative by showing that%
\begin{equation}
1.122n-25<\frac{29+\sqrt{329}}{42}n-25\leq\max_{v\left(  G\right)  =n}\mu
_{1}\left(  G\right)  +\mu_{2}\left(  G\right)  \leq\frac{2}{\sqrt{3}%
}n<1.155n. \label{eq2}%
\end{equation}

\section{Proofs}

Given a graph $G$ and an integer $t>0,$ set $G^{\left[  t\right]  }%
=\overline{\overline{G}^{\left(  t\right)  }},$ i.e., $G^{\left[  t\right]  }$
is obtained from $G^{\left(  t\right)  }$ by joining all vertices within
$V_{u}$ for every $u\in V.$ The following two facts are derived by
straightforward methods.

\emph{(i)} The\ eigenvalues of $G^{\left(  t\right)  }$ are $t\mu_{1}\left(
G\right)  ,\ldots,t\mu_{n}\left(  G\right)  $ together with $n\left(
t-1\right)  $ additional $0$'s.

\emph{(ii)} The\ eigenvalues of $G^{\left[  t\right]  }$ are $t\mu_{1}\left(
G\right)  +t-1,\ldots,t\mu_{n}\left(  G\right)  +t-1$ together with $n\left(
t-1\right)  $ additional $\left(  -1\right)  $'s.

We shall show that the extremal $k$ eigenvalues of $G^{\left(  t\right)  }$
and $G^{\left[  t\right]  }$ are roughly proportional to the corresponding
eigenvalues of $G.$

\begin{lemma}
\label{prop}Let $1\leq k<n,$ $t\geq2.$ Then for every $s\in\left[  k\right]
,$
\begin{align}
0  &  \leq\mu_{s}\left(  G^{\left(  t\right)  }\right)  -t\mu_{s}\left(
G\right)  <\frac{tn}{\sqrt{n-k}},\label{i1}\\
0  &  \geq\mu_{n-s+1}\left(  G^{\left(  t\right)  }\right)  -t\mu
_{n-s+1}\left(  G\right)  >-\frac{tn}{\sqrt{n-k}},\label{i2}\\
0  &  \leq\mu_{s}\left(  G^{\left[  t\right]  }\right)  -t\mu_{s}\left(
G\right)  +t-1<t+\frac{tn}{\sqrt{n-k}},\label{i3}\\
0  &  \geq\mu_{n-s+1}\left(  G^{\left[  t\right]  }\right)  -t\mu
_{n-s+1}\left(  G\right)  +t-1>-t-\frac{tn}{\sqrt{n-k}}. \label{i4}%
\end{align}

\end{lemma}

\begin{proof}
We shall prove (\ref{i1}) first. Fix some $s\in\left[  k\right]  $ and note
that \emph{(i)} implies that $G^{\left(  t\right)  }$ and $G$ have the same
number of positive eigenvalues. In particular, $G^{\left(  t\right)  }$ has at
most $n-1$ negative eigenvalues, and so $\mu_{s}\left(  G^{\left(  t\right)
}\right)  \geq0.$ If $\mu_{s}\left(  G^{\left(  t\right)  }\right)  >0,$ then
$\mu_{s}\left(  G\right)  >0$ and $\mu_{s}\left(  G^{\left(  t\right)
}\right)  =t\mu_{s}\left(  G\right)  ,$ so (\ref{i1}) holds. If $\mu
_{s}\left(  G^{\left(  t\right)  }\right)  =0,$ then
\[
0\geq\mu_{s}\left(  G\right)  \geq\ldots\geq\mu_{n}\left(  G\right)  ,
\]
and inequality (\ref{i1}) follows from
\[
\left(  n-k\right)  \mu_{s}^{2}\left(  G\right)  \leq\left(  n-s+1\right)
\mu_{s}^{2}\left(  G\right)  \leq\sum_{i=s}^{n}\mu_{i}^{2}\left(  G\right)
<n^{2}.
\]

Next we shall prove (\ref{i3}). Note that \emph{(ii)} implies that $G^{\left(
t\right)  }$ and $G$ have the same number of eigenvalues that are greater than
$-1$. Since $G^{\left[  t\right]  }$ has at most $n-1$ eigenvalues that are
less than $-1$, it follows that $\mu_{s}\left(  G^{\left[  t\right]  }\right)
\geq-1.$ If $\mu_{s}\left(  G^{\left[  t\right]  }\right)  >-1,$ then $\mu
_{s}\left(  G\right)  >-1$ and $\mu_{s}\left(  G^{\left[  t\right]  }\right)
=t\mu_{s}\left(  G\right)  +t-1,$ so (\ref{i3}) holds. If $\mu_{s}\left(
G^{\left[  t\right]  }\right)  =-1,$ then
\[
-1\geq\mu_{s}\left(  G\right)  \geq\ldots\geq\mu_{n}\left(  G\right)  ,
\]
and inequality (\ref{i3}) follows from
\[
\left(  n-k\right)  \mu_{s}^{2}\left(  G\right)  <\left(  n-s+1\right)
\mu_{s}^{2}\left(  G\right)  \leq\sum_{i=s}^{n}\mu_{i}^{2}\left(  G\right)
<n^{2}.
\]

Inequalities (\ref{in2}) and (\ref{i4}) follow likewise, with proper changes
of signs.
\end{proof}

We also need the following lemma.

\begin{lemma}
\label{Weyp}Let $G$ be a graph of order $n$ and $H$ be an induced subgraph of
$G$ of order $n-1.$ Then for every $1\leq s\leq3n/4,$%
\begin{align}
0  &  \leq\mu_{s}\left(  G\right)  -\mu_{s}\left(  H\right)  <3\sqrt
{n},\label{Win1}\\
0  &  \geq\mu_{n-s+1}\left(  G\right)  -\mu_{n-s}\left(  H\right)  >-3\sqrt
{n}. \label{Win2}%
\end{align}

\end{lemma}

\begin{proof}
We shall assume that $V\left(  G\right)  =\left\{  1,\ldots,n\right\}  $ and
$V\left(  H\right)  =\left\{  1,\ldots,n-1\right\}  .$ Let $A$ be the
adjacency matrix of $G$ and let $A_{1}$ be the $n\times n$ symmetric matrix
obtained from $A$ by zeroing its $n$th row and column. Since the adjacency
matrix of $H$ is the principal submatrix of $A$ in the first $n-1$ columns and
rows, the\ eigenvalues of $A_{1}$ are $\mu_{1}\left(  H\right)  ,\ldots
,\mu_{n-1}\left(  H\right)  $ together with an additional $0.$ This implies
that, for every $s\in\left[  n-1\right]  ,$%
\begin{equation}
\mu_{s}\left(  A_{1}\right)  =\left\{
\begin{array}
[c]{cc}%
\mu_{s}\left(  H\right)  , & \text{if }\mu_{s}\left(  A_{1}\right)  >0\\
\mu_{s-1}\left(  H\right)  & \text{if }\mu_{s}\left(  A_{1}\right)  \leq0
\end{array}
\right.  \label{obv}%
\end{equation}
We first show that, for every $s\in\left[  n-1\right]  ,$%
\begin{equation}
\mu_{s}\left(  A_{1}\right)  -\mu_{s}\left(  H\right)  \leq\frac{n}{\sqrt
{n-s}}. \label{upb}%
\end{equation}
In view of (\ref{obv}), this is obvious if $\mu_{s}\left(  A_{1}\right)  >0.$
If $\mu_{s}\left(  A_{1}\right)  \leq0,$ again in view of (\ref{obv}), we
have
\[
\mu_{s}\left(  A_{1}\right)  -\mu_{s}\left(  H\right)  =\mu_{s-1}\left(
H\right)  -\mu_{s}\left(  H\right)  \leq\left\vert \mu_{s}\left(  H\right)
\right\vert
\]
Inequality (\ref{upb}) follows now from%
\[
\left(  n-s\right)  \mu_{s}^{2}\left(  G\right)  \leq\left(  n-s+1\right)
\mu_{s}^{2}\left(  G\right)  \leq\sum_{i=s}^{n}\mu_{i}^{2}\left(  G\right)
<n^{2}.
\]
Likewise, with proper changes of signs, we can show that, for every
$s\in\left[  n-1\right]  ,$
\[
\mu_{n-s+1}\left(  A_{1}\right)  -\mu_{n-s}\left(  H\right)  \geq-\frac
{n}{\sqrt{n-s}}.
\]

Having prove (\ref{upb}) we turn to the proof of (\ref{Win1}) and
(\ref{Win2}). Note that the first inequalities in both (\ref{Win1}) and
(\ref{Win2}) follow by Cauchy interlacing theorem. On the other hand, Weyl's
inequalities imply that%
\[
\mu_{n}\left(  A-A_{1}\right)  \leq\mu_{s}\left(  A\right)  -\mu_{s}\left(
A_{1}\right)  \leq\mu_{1}\left(  A-A_{1}\right)  .
\]
Obviously, $\mu_{1}\left(  A-A_{1}\right)  $ is maximal when the off-diagonal
entries of the $n$th row and column of $A$ are $1$'s. Thus, $\mu_{1}\left(
A-A_{1}\right)  \leq\sqrt{n-1}$ and $\mu_{n}\left(  A-A_{1}\right)  =-\mu
_{1}\left(  A-A_{1}\right)  \geq-\sqrt{n-1}.$ Hence,%
\[
\mu_{s}\left(  G\right)  -\mu_{s}\left(  H\right)  =\mu_{s}\left(  A\right)
-\mu_{s}\left(  A_{1}\right)  +\mu_{s}\left(  A_{1}\right)  -\mu_{s}\left(
H\right)  \leq\sqrt{n-1}+\frac{n}{\sqrt{n-s}}<3\sqrt{n}.
\]
Likewise,%
\begin{align*}
\mu_{n-s+1}\left(  G\right)  -\mu_{n-s}\left(  H\right)   &  =\mu
_{n-s+1}\left(  A\right)  -\mu_{n-s+1}\left(  A_{1}\right)  +\mu
_{n-s+1}\left(  A_{1}\right)  -\mu_{n-s}\left(  H\right) \\
&  \geq-\sqrt{n-1}-\frac{n}{\sqrt{n-s}}>-3\sqrt{n},
\end{align*}
completing the proof of Lemma \ref{Weyp}.
\end{proof}

\begin{corollary}
\label{cWeyp}Let $G_{1}$ be a graph of order $n$ and $G_{2}$ be an induced
subgraph of $G_{1}$ of order $n-l.$ Then, for every $1\leq s\leq3\left(
n-l\right)  /4,$%
\begin{align*}
\left\vert \mu_{s}\left(  G_{1}\right)  -\mu_{s}\left(  G_{2}\right)
\right\vert  &  <3l\sqrt{n},\\
\left\vert \mu_{n-s+1}\left(  G_{1}\right)  -\mu_{n-l-s+1}\left(
G_{2}\right)  \right\vert  &  <3l\sqrt{n}.
\end{align*}

\end{corollary}

\begin{proof}
Let $\left\{  v_{1},\ldots,v_{l}\right\}  =V\left(  G_{1}\right)  \backslash
V\left(  G_{2}\right)  .$ Set $H_{0}=G_{1};$ for every $i\in\left[  l\right]
,$ let $H_{i}$ be the subgraph of $G_{1}$ induced by the set $V\left(
G_{1}\right)  \backslash\left\{  v_{1},\ldots,v_{i}\right\}  ;$ clearly,
$H_{l}=G_{2}.$ Since $H_{i+1}$ is an induced subgraph of $H_{i}$ with
$\left\vert H_{i+1}\right\vert =\left\vert H_{i}\right\vert -1,$ Lemma
\ref{Weyp} implies that for every $1\leq s\leq3\left(  n-l\right)  /4,$
\begin{align*}
\left\vert \mu_{s}\left(  G_{1}\right)  -\mu_{s}\left(  G_{2}\right)
\right\vert  &  \leq\sum_{i=0}^{l-1}\left\vert \mu_{s}\left(  H_{i}\right)
-\mu_{s}\left(  H_{i+1}\right)  \right\vert \leq\sum_{i=0}^{l-1}3\sqrt
{n-i}<3l\sqrt{n},\\
\left\vert \mu_{n-s+1}\left(  G_{1}\right)  -\mu_{n-l-s+1}\left(
G_{2}\right)  \right\vert  &  \leq\sum_{i=0}^{l-1}\left\vert \mu
_{n-i+s+1}\left(  H_{i}\right)  -\mu_{n-i-1-s+1}\left(  H_{i+1}\right)
\right\vert \\
&  \leq\sum_{i=0}^{l-1}3\sqrt{n-i}<3l\sqrt{n},
\end{align*}
completing the proof of the corollary.
\end{proof}

\begin{proof}
[\textbf{Proof of Theorem \ref{th1} }]Set
\[
\varphi\left(  n\right)  =\frac{1}{n}\max\left\{  F\left(  G\right)
:G\in\mathcal{F}\text{, }v\left(  G\right)  =n\right\}
\]
Let $M=\sum_{i=1}^{k}\left\vert \alpha_{i}\right\vert +\left\vert \beta
_{i}\right\vert +\left\vert \gamma_{i}\right\vert +\left\vert \delta
_{i}\right\vert $ and set
\[
c=\lim_{n\rightarrow\infty}\sup\varphi\left(  n\right)  .
\]
Since $\left\vert F\left(  G\right)  \right\vert \leq Mn,$ the value $c$ is
defined. We shall prove that, in fact, $c$ satisfies (\ref{eq1}).

Note first if $t\geq2$, and $n>4k/3,$ then, for any $i\in\left[  k\right]  ,$
Lemma \ref{prop} implies that
\begin{equation}
F\left(  G^{\left(  t\right)  }\right)  -tF\left(  G\right)  \geq-M\left(
t+\frac{tn}{\sqrt{n-k}}\right)  \geq-M\left(  t+2t\sqrt{n}\right)
\geq-3Mt\sqrt{n}. \label{in1}%
\end{equation}

Select $\varepsilon>0$ and let $G\in\mathcal{F}$ be a graph of order
$n>\left(  3M/\varepsilon\right)  ^{2}$ such that
\[
c+\varepsilon\geq\varphi\left(  n\right)  =\frac{F\left(  G\right)  }{n}\geq
c-\varepsilon.
\]
Suppose $N\geq n\left\lceil n\max\left\{  2,\left(  \left\vert c\right\vert
/\varepsilon+1\right)  ,\left(  3M/\varepsilon\right)  ^{2}\right\}
\right\rceil ;$ therefore the value $t=\left\lfloor N/n\right\rfloor $
satifies $t\geq n\max\left\{  2,\left(  \left\vert c\right\vert /\varepsilon
+1\right)  ,\left(  3M/\varepsilon\right)  ^{2}\right\}  .$ We shall show that
$\varphi\left(  N\right)  \geq c-4\varepsilon,$ which implies the assertion.

Let $G_{1}$ be the union of $G^{\left(  t\right)  }$ and $N-tn$ isolated
vertices. Clearly $v\left(  G_{1}\right)  =N$ and, since $\mathcal{F}$ is
multiplicative, $G_{1}\in\mathcal{F}$. In view of $N-tn<n,$ Corollary
\ref{cWeyp} implies that%
\[
F\left(  G_{1}\right)  \geq F\left(  G^{\left(  t\right)  }\right)
-3Mn\sqrt{N}.
\]
Therefore, in view of $\varphi\left(  N\right)  \geq F\left(  G_{1}\right)
/N$ and (\ref{in1}),
\[
\varphi\left(  N\right)  \geq\frac{F\left(  G^{\left(  t\right)  }\right)
-3Mn\sqrt{N}}{N}\geq\frac{tF\left(  G\right)  -3Mt\sqrt{n}-3Mn\sqrt{N}}{N}.
\]
We find that
\begin{align*}
\varphi\left(  N\right)   &  \geq\frac{F\left(  G\right)  }{n}-\frac{F\left(
G\right)  \left(  N-tn\right)  }{nN}-\frac{3Mt\sqrt{n}+3Mn\sqrt{N}}{N}\\
&  \geq\varphi\left(  n\right)  -\frac{\left\vert \varphi\left(  n\right)
\right\vert n}{N}-\frac{3Mt\sqrt{n}+3Mn\sqrt{N}}{N}\\
&  \geq\varphi\left(  n\right)  -\frac{n^{2}\left(  \left\vert c\right\vert
+\left\vert \varepsilon\right\vert \right)  +3Mt\sqrt{n}+3Mn\sqrt{N}}{N}\\
&  \geq\varphi\left(  n\right)  -\frac{n^{2}\left(  \left\vert c\right\vert
+\left\vert \varepsilon\right\vert \right)  +3Mt\sqrt{n}}{nt}-\frac{3Mn}%
{\sqrt{nt}}\\
&  =\varphi\left(  n\right)  -\frac{n\left(  \left\vert c\right\vert
+\left\vert \varepsilon\right\vert \right)  }{t}-\frac{3M}{\sqrt{n}}%
-3M\sqrt{\frac{n}{t}}\geq c-4\varepsilon,
\end{align*}
completing the proof of Theorem \ref{th1}.
\end{proof}

We turn now to the proof of inequality (\ref{eq2}); we present it in two propositions.

\begin{proposition}
If $G$ is a graph of order $n$, then $\mu_{1}\left(  G\right)  +\mu_{2}\left(
G\right)  \leq\left(  2/\sqrt{3}\right)  n.$
\end{proposition}

\begin{proof}
Setting $m=e\left(  G\right)  ,$ we see that
\begin{equation}
\mu_{1}^{2}\left(  G\right)  +\mu_{2}^{2}\left(  G\right)  \leq\mu_{1}%
^{2}\left(  G\right)  +\ldots+\mu_{n}^{2}\left(  G\right)  =2m.\label{in2}%
\end{equation}
If $m\leq n^{2}/4,$ the result follows from
\[
\mu_{1}\left(  G\right)  +\mu_{2}\left(  G\right)  \leq\sqrt{2\left(  \mu
_{1}^{2}\left(  G\right)  +\mu_{2}^{2}\left(  G\right)  \right)  }\leq
2\sqrt{m}\leq n,
\]
so we shall assume that $m>n^{2}/4.$ From (\ref{in2}), we clearly have
\[
\mu_{1}\left(  G\right)  +\mu_{2}\left(  G\right)  \leq\sqrt{2m-\mu_{2}%
^{2}\left(  G\right)  }+\mu_{2}\left(  G\right)  .
\]
The value $\sqrt{2m-x^{2}}+x$ is increasing in $x$ for $x\leq m.$ On the other
hand, Weyl's inequalities imply that
\[
\mu_{2}\left(  G\right)  +\mu_{n}\left(  \overline{G}\right)  \leq\mu
_{2}\left(  K_{n}\right)  =-1.
\]
Hence, if $G\neq K_{n},$ we have $\mu_{2}\left(  G\right)  \geq0$ and so,
$\mu_{2}^{2}\left(  G\right)  <\mu_{n}^{2}\left(  \overline{G}\right)  ;$ if
$G=K_{n},$ then $\mu_{2}^{2}\left(  G\right)  =\mu_{n}^{2}\left(  \overline
{G}\right)  +1;$ thus we always have
\[
\mu_{2}^{2}\left(  G\right)  \leq\mu_{n}^{2}\left(  \overline{G}\right)  +1.
\]
From
\[
\mu_{2}^{2}\left(  G\right)  \leq\mu_{n}^{2}\left(  \overline{G}\right)
+1\leq e\left(  \overline{G}\right)  +1\leq\frac{n\left(  n-1\right)  }%
{2}+1-m\leq\frac{n^{2}}{2}-m<m,
\]
we see that
\[
\mu_{1}\left(  G\right)  +\mu_{2}\left(  G\right)  \leq\sqrt{3m-n^{2}/2}%
+\sqrt{n^{2}/2-m}.
\]
The right-hand side of this inequality is maximal for $m=n^{2}/3$ and the
result follows.
\end{proof}

\begin{proposition}
For every $n\geq21$ there exists a graph of order $n$ with
\[
\mu_{1}\left(  G\right)  +\mu_{2}\left(  G\right)  >\frac{29+\sqrt{329}}%
{42}n-25.
\]

\end{proposition}

\begin{proof}
Suppose $n\geq21,$ set $k=\left\lfloor n/21\right\rfloor ,$ let $G_{1}$ be the
union of two copies of $K_{8k}$ and $G_{2}$ be the join of $K_{5k}$ and
$G_{1};$ clearly $v\left(  G_{2}\right)  =21k.$ Add $n-21k$ isolated vertices
to $G_{2}$ and write $G$ for the resulting graph. By Cauchy interlacing
theorem, we have
\begin{align*}
\mu_{1}\left(  G\right)   &  \geq\mu_{1}\left(  G_{2}\right)  ,\\
\mu_{2}\left(  G\right)   &  \geq\mu_{2}\left(  G_{2}\right)  \geq\mu
_{2}\left(  G_{1}\right)  =8k-1.
\end{align*}
Since the graphs $K_{5k}$ and $G_{1}$ are regular, a theorem of Finck and
Grohmann \cite{FiGr65} (see also \cite{CDS80}, Theorem 2.8) implies that
$\mu_{1}\left(  G_{2}\right)  $ is the positive root of the equation
\[
\left(  x-5k+1\right)  \left(  x-8k+1\right)  -80k^{2}=0.
\]
Hence,
\begin{align*}
\mu_{1}\left(  G\right)  +\mu_{2}\left(  G\right)   &  \geq\frac{\left(
29k-4\right)  +k\sqrt{329}}{2}>\frac{\left(  29\left(  n-20\right)
-84\right)  +\left(  n-20\right)  \sqrt{329}}{42}\\
&  >\frac{29+\sqrt{329}}{42}n-25,
\end{align*}
completing the proof.
\end{proof}

\textbf{Acknowledgments }Part of this research was completed while the author
was visiting the Institute for Mathematical Sciences, National University of
Singapore in 2006. The author is also indebted to B\'{e}la Bollob\'{a}s for
his kind support.

Finally, the referees' criticisms helped to correct the first version of the paper.

\end{document}